\documentclass[11pt]{amsart}
\input rlepsf

\newcommand{\U}{{\mbox{\rm U}}}
\newcommand{\SO}{{\mbox{\rm SO}}}
\newcommand{\OO}{{\mbox{\rm O}}}
\newcommand{\ph}{{\mbox{\rm ph}}}
\newcommand{\bfz}{{\mathbb {Z}}}
\newcommand{\bfc}{{\mathbb {C}}}

\newcommand{\num}{\mbox{\rm numerator}}
\newcommand{\den}{\mbox{\rm denominator}}
\newcommand{\fo}{{\mathfrak o}}

\begin{document}

\title{$E_8$-plumbings and exotic contact structures on spheres}
\author{Fan Ding}
\author{Hansj\"org Geiges}
\address{Department of Mathematics, Peking University, Beijing 100871, 
P.R. China}
\email{dingfan@math.pku.edu.cn}
\address{Mathematisches Institut, Universit\"at zu K\"oln, 
Weyertal 86-90, 50931 K\"oln, Germany}
\email{geiges@math.uni-koeln.de}
\date{}
\maketitle

\newtheorem{thm}{Theorem}
\newtheorem{lem}[thm]{Lemma}
\newtheorem{prop}[thm]{Proposition}

\theoremstyle{definition}
\newtheorem*{rem}{Remark}
\newtheorem*{rems}{Remarks}
\newtheorem*{ack}{Acknowledgements}

\section{Introduction}
The standard contact structure $\xi_{st}$ on the unit
sphere $S^{2n-1}=\partial D^{2n}\subset {\mathbb C}^n$
can be defined as the hyperplane field of complex tangencies. In other
words, if we write $J_0$ for the complex structure on (the tangent bundle of)
${\mathbb C}^n$, then $\xi_{st}(p)=T_pS^{2n-1}\cap J_0(T_pS^{2n-1})$ for
all $p\in S^{2n-1}$.

Conversely, given any contact structure $\xi =\ker\alpha$ on $S^{2n-1}$
(with $\alpha$ a $1$-form such that $\alpha\wedge (d\alpha )^{n-1}$
is a volume form defining the standard orientation),
there is a homotopically unique complex bundle structure $J$ on $\xi$
such that $d\alpha (.,J.)$ is a $J$-invariant Riemannian metric on $\xi$.
This extends to a complex structure $J$ on $T{\mathbb R}^{2n}|_{S^{2n-1}}$
by requiring that $J$ send the outer normal of $S^{2n-1}$ to a vector field
$X$ on $S^{2n-1}$ with $\alpha (X)>0$. If this $J$ extends
as an almost complex structure over the disc $D^{2n}$, then $\xi$
is called {\it homotopically trivial}.

The complex structure $J|_{\xi}$ and the trivial line bundle
spanned by the vector field $X$ can be interpreted as a reduction
of the structure group of $TS^{2n-1}$ to $\U_{n-1}\times 1$. Such a reduction
is called an {\it almost contact structure}. Homotopy classes of
almost contact structures on $S^{2n-1}$ are classified by
\[
\pi_{2n-1}(\SO_{2n-1}/\U_{n-1})\cong
\left\{ \begin{array}{ll}
\bfz\oplus\bfz_2&\;\; \mbox{\rm for}\; n\equiv 0\;\mbox{\rm mod}\; 4,\\
\bfz_{(n-1)!}   &\;\; \mbox{\rm for}\; n\equiv 1\;\mbox{\rm mod}\; 4,\\
\bfz            &\;\; \mbox{\rm for}\; n\equiv 2\;\mbox{\rm mod}\; 4,\\
\bfz_{(n-1)!/2} &\;\; \mbox{\rm for}\; n\equiv 3\;\mbox{\rm mod}\; 4.
\end{array}\right.
\]

Homotopically trivial
contact structures are those that induce the same almost contact structure
as $\xi_{st}$ and correspond to the trivial element in the respective
group above.

The purpose of this note is to prove the following result, where we
call a contact structure on $S^{2n-1}$ {\it exotic} if it is not
diffeomorphic to $\xi_{st}$.

\begin{thm}
\label{thm:contact}
In any odd dimension $2n-1\geq 3$,
the standard sphere $S^{2n-1}$ admits exotic but homotopically
trivial contact structures.
\end{thm}

For the $3$-sphere this was proved by Bennequin~\cite{benn83}.
The (by \cite{elia89} essentially unique) exotic but homotopically
trivial contact structure on $S^3$ is not symplectically fillable.

In higher dimensions, by contrast, exoticity of $(S^{2n-1},\xi )$
can be shown by exhibiting a symplectic filling
of $(S^{2n-1},\xi )$ that does not contain symplectic $2$-spheres
and is not diffeomorphic to a disc~$D^{2n}$. This result is due
to Eliashberg, Floer, Gromov and McDuff, see~\cite{elia91}.

In that paper, Eliashberg used this result to derive Theorem~\ref{thm:contact}
for spheres of dimension $4k+1$.
Notice that this corresponds to $n=2k+1$ odd,
in which case there are only finitely many homotopy classes
of almost contact structures. For that reason, it is enough to find
{\em some} exotic $(S^{2n-1},\xi )$ (with a symplectic
filling as described) --- Eliashberg did this by using a plumbing
construction. An exotic but homotopically
trivial contact structure on $S^{2n-1}$
can then be constructed simply by taking the
connected sum of suitably many copies of $(S^{2n-1},\xi )$: the result
of \cite{elia91} applies to the boundary connected sum of the fillings.

For $n$ even, this approach fails. In~\cite{geig97} a more careful study
of the homotopy classes of almost contact structures was undertaken, and
Theorem~\ref{thm:contact} was proved for $S^7$ and spheres of dimension
$\equiv 3$ mod~$8$. The general case of spheres of dimension $\equiv 7$
mod~$8$ (that is, $n\equiv 0$ mod~$4$), however, had remained elusive.
In the present note, we show that these remaining cases can be settled
by extending ideas from~\cite{elia91} and~\cite{geig97}.
Starting from an $E_8$-plumbing, we show how to keep control of the
homotopical information necessary to prove Theorem~\ref{thm:contact}.
In the process, we construct almost complex manifolds that are
of some independent interest.

The `classical' methods of the present note do not allow us to
distinguish different exotic contact structures in the same homotopy
class of almost contact structures. For that, contact homological methods
are required, see \cite{usti99}, \cite{bour03} (again, those papers
only deal with spheres of dimension~$4k+1$).

For general homotopy classes we have the following statement.

\begin{thm}
\label{thm:contact2}
The standard spheres of dimension $\equiv 1,3,5$ {\rm mod}~$8$ (and greater
than~$1$) admit exotic contact structures in every homotopy class
of almost contact structures; so does~$S^7$. The spheres of dimension
$\equiv 7$ {\rm mod}~$8$ admit exotic contact structures in every
stably trivial homotopy class of almost contact structures.
\end{thm}

For dimension $3$ this is due to Lutz~\cite{lutz70}; exoticity follows from
these structures being overtwisted. For dimensions $2n-1\equiv 1,5$ mod~$8$
(and implicitly for $2n-1=7$ or $\equiv 3$ mod~$8$) this was proved
in~\cite{geig97}.
The exoticity of the structures follows from the same fillability
argument as above. (A contact structure in a nontrivial homotopy class
of almost contact structures need not be exotic, {\it a priori}:
the diffeomorphism group of $S^{2n-1}$ may act nontrivially on these
homotopy classes.)

We do not know anything about the realisation of stably nontrivial
homotopy classes of almost contact structures on spheres of dimension
$\equiv 7$ mod~$8$ (and greater than~$7$). At the end of
this note we show that, at least on $S^{15}$, such homotopy classes
cannot be realised by contact structures that are symplectically
fillable by some highly connected manifold. In other words, the
construction used to prove Theorems \ref{thm:contact}
and~\ref{thm:contact2} does not allow one to extend
Theorem~\ref{thm:contact2} to all homotopy classes of almost
contact structures.
\section{The $E_8$-plumbing}
Let $DTS^{2m}$ be the tangent unit disc bundle of the $2m$-dimensional sphere,
$m\geq 2$. It is well-known that the plumbing of eight copies of
$DTS^{2m}$ according to the $E_8$-graph

\begin{figure}[h]
\centerline{\relabelbox\small
\epsfxsize 8cm \epsfbox{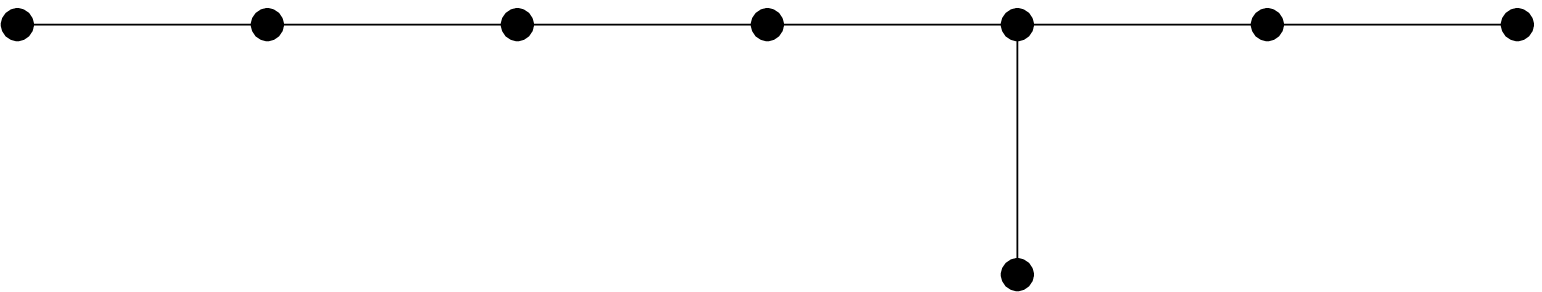}
\endrelabelbox}
\end{figure}

\noindent
yields a manifold $V^{4m}$ with boundary $\partial V^{4m}=\Sigma^{4m-1}$
a homotopy sphere, cf.\ \cite{brow72}, \cite[\S 8]{hima68}. The intersection
matrix of $V^{4m}$ (with the natural orientation of that manifold)
is the $E_8$-matrix
\[\left(\begin{array}{cccccccc}
2 & 1 &   &   &   &   &   &   \\
1 & 2 & 1 &   &   &   &   &   \\
  & 1 & 2 & 1 &   &   &   &   \\
  &   & 1 & 2 & 1 &   &   &   \\
  &   &   & 1 & 2 & 1 & 0 & 1 \\
  &   &   &   & 1 & 2 & 1 & 0 \\
  &   &   &   & 0 & 1 & 2 & 0 \\
  &   &   &   & 1 & 0 & 0 & 2 
\end{array}\right)\]
(with zeros in the blank spaces). In fact, $\Sigma^{4m-1}$ is a generator
of the group $bP_{4m}$ of homotopy $(4m-1)$-spheres bounding a parallelisable
manifold, see \cite{kemi63} and \cite[\S 10]{hima68}. The order
$|bP_{4m}|$ of this group is given by
\[ |bP_{4m}|=2^{2m-2}(2^{2m-1}-1)\cdot\num\left(\frac{4B_m}{m}\right) ,\]
where $B_m$ denotes the $m$th Bernoulli number. (In \cite{kemi63}
this was only proved up to a potential factor of $2$ in the case
$m$ even; this was later settled by Mahowald~\cite{maho70},
cf.~\cite[p.~771]{hirz87}.) The denominator of $B_m$ is square-free
and divisible by~$2$ \cite[p.~284]{mist74}. This implies
\[ |bP_{4m}|=2^{2m-2}(2^{2m-1}-1)\cdot a_m\cdot \num\left(\frac{B_m}{4m}\right) ,\]
where $a_m=2$ for $m$ odd, $a_m=1$ for $m$ even.

The boundary connected sum of $|bP_{4m}|$ copies of $V^{4m}$ has
boundary $S^{4m-1}$. Write $W_0^{4m}$ for the closed, orientable $4m$-manifold
obtained by attaching a $4m$-disc $D^{4m}$ along the boundary.
(There is a choice here, since for some $m$ there are diffeomorphisms
of $S^{4m-1}$ that do not extend over $D^{4m}$, but the resulting
ambiguity is of no consequence for our considerations.)

We collect some well-known topological information about $W_0^{4m}$,
cf.~\cite[p.~770]{hirz87}. The manifold $W_0^{4m}$
is $(2m-1)$-connected and its middle cohomology group is
$H^{2m}(W_0^{4m})\cong \bfz^{\oplus 8|bP_{4m}|}$. The plumbing $V^{4m}$
has as deformation retract a union of eight $2m$-spheres intersecting one
another according to the $E_8$-graph. The restriction of the tangent
bundle $TV^{4m}$ to each of these spheres is isomorphic to the
Whitney sum $TS^{2m}\oplus TS^{2m}$, which is stably trivial. It follows
that $TV^{4m}$ is trivial, and so is $T(W_0^{4m}-D^{4m})$. Hence, for
$m$ even, the Pontrjagin class $p_{m/2}(W_0^{4m})$ is trivial.

\begin{lem}
\label{lem:pont}
The Pontrjagin number $p_m[W_0^{4m}]$ is given by
\[ p_m[W_0^{4m}] = a_m\cdot \den\left(\frac{B_m}{4m}\right)\cdot (2m-1)!,\]
where $a_m=2$ for $m$ odd, $a_m=1$ for $m$ even.
\end{lem}

\begin{proof}
The signature of the $E_8$-matrix is equal to~$8$, so the signature
of $W_0^{4m}$ equals
\begin{equation}
\label{eqn:s1}
\sigma (W_0^{4m})=8|bP_{4m}|=2^{2m+1}(2^{2m-1}-1)\cdot a_m\cdot \num
\left(\frac{B_m}{4m}\right).
\end{equation}
The coefficient of $p_m$ in the $L$-polynomial $L_m(p_1,\ldots ,p_m)$
is equal to
\[ d_m:=2^{2m}(2^{2m-1}-1)\frac{B_m}{(2m)!},\]
see \cite[p.~12]{hirz78}. By the signature theorem of
Hirzebruch~\cite[Thm.~8.2.2]{hirz78},
\begin{equation}
\label{eqn:s2}
\sigma (W_0^{4m})=d_mp_m[W_0^{4m}].
\end{equation}
The formula for $p_m[W_0^{4m}]$ follows by comparing the expressions
(\ref{eqn:s1}) and (\ref{eqn:s2}) for $\sigma (W_0^{4m})$.
\end{proof}

\section{Almost complex manifolds} 
Starting from the manifold $W_0^{4m}$ of the preceding section, we now want to
construct almost complex $4m$-manifolds with certain special features.
To do this, we shall need to analyse the top-dimensional obstruction to
the existence of almost complex structures. Strictly speaking, only
the case $m$ even will be relevant for the application to contact
geometry, but the construction for $m$ odd is analogous, and of interest
because of the following result.

\begin{prop}
\label{prop:acs}
The Euler characteristic of an almost complex manifold of (real) dimension
$4m$ with first Chern class $c_1\equiv 0$ {\rm mod}~$2$ and vanishing
decomposable Chern numbers is divisible by
\[ a_m\cdot\den\left(\frac{B_m}{2m}\right)\cdot (2m-1)!.\]
There are almost complex manifolds $M_0^{4m}$ with the described
properties that realise this minimal (non-zero) absolute value of
the Euler characteristic.
\end{prop}

\begin{rems}
(1) Such manifolds have been found previously by Puschnigg (unpublished),
cf.~\cite[p.~778]{hirz87}; his diploma thesis~\cite{pusc88}
contains a homotopy-theoretic proof of (essentially) the above
proposition, but no explicit construction of the~$M_0^{4m}$.

(2) If only the vanishing of the decomposable Chern numbers is
required, the formula for the minimal positive Euler characteristic
is as above, but without the factor~$a_m$, see~\cite{brum68}
and~\cite{bake83}. No explicit manifolds realising this value were given
in those papers.
\end{rems}

\begin{proof}
The divisibility statement is contained in~\cite[Lemma~2.3]{brum68}.
Since the proof there is left to the reader, here is a brief
sketch, which is intended merely to indicate what integrality statements
enter into the proof. We assume for the time being that the reader
is familiar with the notations for various classical genera of manifolds;
in the more formal arguments below all these notations will be
explained.

Let $M$ be a manifold satisfying the assumptions of the proposition.
The condition on $c_1$ implies that $M$ is a spin manifold, hence
$\hat{A}[M]\in a_m\bfz$, see~\cite[p.~198/9]{hirz78}. Because of the
vanishing of the decomposable Chern numbers, only the coefficient of
the Pontrjagin class $p_m$ in the $\hat{A}$-polynomial is
relevant for the computation of the genus $\hat{A}[M]$. The
way to compute that coefficient is described in
\cite[Section~1.4]{hirz78}, and one finds
\[ \hat{A}[M]=-\frac{1}{(2m-1)!}\cdot\frac{B_m}{4m}\cdot p_m[M].\]
We conclude that
\[ a_m\cdot (2m-1)!\cdot\den\left(\frac{B_m}{4m}\right)\; |\;
\num\left(\frac{B_m}{4m}\right)\cdot p_m[M].\]
From the integrality of $\hat{A}[M]$ and
$\langle \ph (M)\hat{A}(M),[M]\rangle$ for the spin manifold~$M$, and the
vanishing of the decomposable Chern numbers, one finds with a little
computation that
\[ (2m-1)!\; |\; p_m[M].\]
Hence
\[ a_m\cdot (2m-1)!\cdot\den\left(\frac{B_m}{4m}\right)\; |\;
p_m[M].\]
The factor $a_m$ in that last divisibility statement is justified
by the fact that the numerator of $B_m/4m$ is odd; we do not
need to use the stronger statement that for the spin manifold $M$ the genus
$\langle \ph (M)\hat{A}(M),[M]\rangle$ is an integer divisible by~$a_m$.

Finally, the vanishing of the decomposable Chern numbers implies
that $p_m[M]$ is (up to sign)
twice the Euler characteristic of~$M$. This yields the
claimed divisibility of the Euler characteristic. 

We now turn to the construction of the manifolds~$M_0^{4m}$.
Let $W$ be an oriented $4m$-dimensional manifold, and assume that there
exists an almost complex structure $J$
on $W-D^{4m}$ for some embedded disc~$D^{4m}$. Write
\[ \fo (W,J)\in H^{4m}(W;\pi_{4m-1}(\SO_{4m}/\U_{2m}))
\cong\left\{\begin{array}{cl}
\bfz\oplus\bfz_2 & \;\;\mbox{\rm for}\; m\equiv 0\;\mbox{\rm mod}\; 2,\\
\bfz             & \;\;\mbox{\rm for}\; m\equiv 1\;\mbox{\rm mod}\; 2,
\end{array}\right.\]
for the obstruction to extending $J$ as an almost complex structure
over~$W$. Here the splitting $\pi_{8k-1}(\SO_{8k}/U_{4k})\cong\bfz
\oplus\bfz_2$ is defined by identifying the $\bfz$-summand with the
kernel of the stabilising map
\[ \pi_{8k-1}(\SO_{8k}/U_{4k})\longrightarrow\pi_{8k-1}(\SO/\U )\cong
\bfz_2,\]
cf.~\cite[Lemma~8]{kahn69}.
For $m\equiv 0$ mod~$2$ we write $\fo =\fo_0+\fo_2\in\bfz\oplus\bfz_2$;
for $m\equiv 1$ mod~$2$ we also write $\fo =\fo_0$ in order to simplify
notation below.

Then we have the following statements from~\cite{kahn69} and~\cite{mass61}:
\begin{itemize}
\item[(i)] $\fo (S^{4m},J)$ is independent of $J$ and will be written
as $\fo (S^{4m})$. We have $\fo (S^{8k})=(1,0)$ and $\fo (S^{8k-4})=1$.
\item[(ii)] Almost complex structures $J$ on $W-D^{4m}$ and $J'$ on
$W'-D^{4m}$ give rise to a natural almost complex structure $J+J'$
on $W\# W'-D^{4m}$ (which coincides with $J$ or $J'$ along the
$(4m-1)$-skeleton of $W$ or $W'$, respectively) such that
\[ \fo (W\# W',J+J')=\fo (W,J)+\fo (W',J')-\fo (S^{4m}).\]
\item[(iii)] Write $p_m$ for the top-dimensional Pontrjagin class of the
tangent bundle of $W$ and $c_i$ for the Chern classes of the complex
bundle $(T(W-D^{4m}),J)$. Then the obstruction $\fo_0 (W,J)$ can be computed
explicitly as
\[ \fo_0(W,J)=\frac{1}{2}\chi (W)+\frac{1}{4}\langle (-1)^{m+1}p_m+
\sum_{\substack{i+j=2m\\ i,j\geq 1}} (-1)^ic_ic_j,[W]\rangle ,\]
where $\chi (W)$ denotes the Euler characteristic of~$W$.
If $J$ extends over $W$ as a stable almost complex structure $\widetilde{J}$
with top-dimensional Chern class $c_{2m}(\widetilde{J})$, this
formula simplifies to
\[ \fo_0(W,J)=\frac{1}{2} \left( \chi (W)-
\langle c_{2m}(\widetilde{J}),[W]\rangle \right) .\]
\end{itemize}

Since the tangent bundle $T(W_0^{4m}-D^{4m})$ is trivial, the
manifold $W_0^{4m}-D^{4m}$ admits an almost complex structure $J_0$
inducing the natural orientation and with vanishing Chern
class~$c_m$. Thus, we obtain
\begin{eqnarray*}
\fo_0(W_0^{4m},J_0) & = & \frac{1}{2}\,\chi (W_0^{4m})+
\frac{(-1)^{m+1}}{4}\, p_m[W_0^{4m}]\\
  & = & 1+\frac{1}{2}\left( 1+\frac{(-1)^{m+1}}{2d_m}\right)\sigma (W_0^{4m}).
\end{eqnarray*}

(a) For $m=2k-1$ an odd integer, $q:=\fo (W_0^{8k-4},J_0)$ is a positive
integer. The parallelisable manifold $S^1\times S^{8k-5}$ admits an
almost complex structure $J'$, hence $\fo (S^1\times S^{8k-5},J')=0$.
Define
\[ M_0^{8k-4}=W_0^{8k-4}\#_q(S^1\times S^{8k-5}).\]
By (i) and (ii), the manifold $M_0^{8k-4}$ admits an almost complex structure.
The Euler characteristic of this manifold is given by
\begin{eqnarray*}
\chi (M_0^{8k-4}) & = & 2+\sigma (W_0^{8k-4})-2q\\
                  & = & -\frac{1}{2}\cdot p_{2k-1}[W_0^{8k-4}]\\
                  & = & -2\cdot \den\left( \frac{B_{2k-1}}{2(2k-1)}\right)
                        \cdot (2(2k-1)-1)!. 
\end{eqnarray*}

(b) For $m=2k$ an even integer, we want to show that $q:=\fo_0(W_0^{8k},J_0)$
is a negative integer. A famous formula of Euler, cf.~\cite[p.~286]{mist74},
states that
\[ \frac{B_m(2\pi )^{2m}}{2(2m)!}=\sum_{\nu=1}^{\infty}\frac{1}{\nu^{2m}}.\]
Hence, with $d_m$ as defined in the proof of Lemma~\ref{lem:pont},
we get (with $m\geq 2$)
\begin{eqnarray*}
d_m & =    & \frac{2(2^{2m-1}-1)}{\pi^{2m}}\, \sum_{\nu =1}^{\infty}
                      \frac{1}{\nu^{2m}}\\
    & \leq & \frac{2(2^{2m-1}-1)}{\pi^{2m}}\, \sum_{\nu =1}^{\infty}
                      \frac{1}{\nu^2}\\
    & =    & \frac{2(2^{2m-1}-1)}{\pi^{2m}}\cdot\frac{\pi^2}{6}\\
    & \leq & \bigl(\frac{2}{\pi}\bigr)^{2m}\cdot 2
             \; \leq \;
             \bigl(\frac{2}{3}\bigr)^4\cdot 2
             \; = \;
             \frac{32}{81}.
\end{eqnarray*}
It follows that
\[ \frac{1}{2}-\frac{1}{4d_{2k}}\leq\frac{1}{2}-\frac{81}{128}=
 -\frac{17}{128}\]
and
\[ q=1+\bigl(\frac{1}{2}-\frac{1}{4d_{2k}}\bigr)\sigma (W_0^{8k})
   \leq 1-\frac{17}{128}\cdot 8<0.\]

Define
\[ M_0^{8k}=W_0^{8k}\#_{-q}(S^{4k}\times S^{4k}).\]

\begin{lem}
\label{lem:M_0}
The manifold $M_0^{8k}$ admits an almost complex structure.
\end{lem}

\begin{proof}
Since the manifold $S^{4k}\times S^{4k}$ has stably trivial tangent
bundle, we can find an almost complex structure $J'$ on
$S^{4k}\times S^{4k}-D^{8k}$ that extends as a stable structure
with vanishing Chern classes. Hence $\fo_0(S^{4k}\times S^{4k},J')=2$
by (iii) and $\fo_0(M_0^{8k},J_0+(-q)J')=0$ by (i) and (ii).

The vanishing of $\fo_2$ can be deduced from~\cite[Thm.~1.2]{heap70}.
This theorem states that subject to the condition $\mbox{\rm Sq}^2
H^{8k-2}(M_0^{8k})=0$, which is obviously satisfied, we have
$\fo_2=0$ if and only if
\begin{eqnarray}
\label{eqn:heaps}
\langle \ph (M_0^{8k})\hat{A}(M_0^{8k}),[M_0^{8k}]\rangle\equiv 0
\;\mbox{\rm mod}\; 2.
\end{eqnarray}
Here $\ph (M_0^{8k})$ denotes the Pontrjagin character of $TM_0^{8k}$,
i.e.\ the Chern character of $TM_0^{8k}\otimes\bfc$, and
$\hat{A}(M_0^{8k})=\sum_{\nu =0}^{\infty}\hat{A}_{\nu}(p_1,\ldots ,p_{\nu})$
the full $\hat{A}$-polynomial in the Pontrjagin classes of~$M_0^{8k}$. Write
\[ \ph (M_0^{8k})=8k+\frac{s_{2k}}{(2k)!}+\frac{s_{4k}}{(4k)!}\]
with $s_j\in H^{2j}(M_0^{8k})$, cf.~\cite[p.~92]{hirz78}. From the
Newton formulae
\[ s_j-c_1s_{j-1}+\ldots +(-1)^j c_j\cdot j=0,\]
where the $c_i$ are the Chern classes of $TM_0^{8k}\otimes\bfc$, we find
$s_{2k}=-2k\cdot c_{2k}=(-1)^{k+1}2k\cdot p_k(M_0^{8k})=0$
and $s_{4k}=-4k\cdot c_{4k}=-4k\cdot p_{2k}(M_0^{8k})$. Hence
\[ \ph (M_0^{8k})=8k-\frac{p_{2k}(M_0^{8k})}{(4k-1)!},\]
and therefore
\[ \langle \ph (M_0^{8k})\hat{A}(M_0^{8k}),[M_0^{8k}]\rangle =
8k\cdot \hat{A}[M_0^{8k}]-\frac{p_{2k}[M_0^{8k}]}{(4k-1)!}.\]
The first summand is even, since $\hat{A}[M_0^{8k}]$ is an integer.
By the proof of Lemma~\ref{lem:pont} --- notice that $\sigma (M_0^{8k})=
\sigma (W_0^{8k})$ ---, the second summand equals the denominator of
$B_{2k}/8k$, which is even, and thus indeed $\fo_2=0$ by the
cited result of~\cite{heap70}.
\end{proof}

\begin{rem}
Concerning the vanishing of the stable obstruction~$\fo_2$,
Pusch\-nigg~\cite{pusc04} has shown us a more general $K$-theoretic argument
for the existence of a stable almost complex structure on every
almost parallelisable manifold of even (real) dimension.
\end{rem}

The Euler characteristic of $M_0^{8k}$ equals
\begin{eqnarray*}
\chi (M_0^{8k}) & = & 2+\sigma (W_0^{8k})-2q\\
                & = & \frac{1}{2}\cdot p_{2k}[W_0^{8k}]\\
                & = & \den\left(\frac{B_{2k}}{2\cdot 2k}\right)
                        \cdot (2\cdot 2k-1)!.
\end{eqnarray*}
This concludes the proof of Proposition~\ref{prop:acs}.
\end{proof}

\section{Proof of Theorems \ref{thm:contact} and \ref{thm:contact2}}
Theorem~\ref{thm:contact} for spheres of dimension $\equiv 7$ mod~$8$
follows from the existence of the $(4k-1)$-connected almost
complex $8k$-manifold $M_0^{8k}$ via \cite[Thm.~24]{geig97}. Given what
was known previously, this completes the proof of Theorem~\ref{thm:contact}.

Theorem~\ref{thm:contact2} follows similarly. Write $J$ for the almost
complex structure on $M_0^{8k}$. Then obviously $\fo (M_0^{8k},J)=0$.
By ~\cite{kahn69}, the manifold $-M_0^{8k}$ (that is, $M_0^{8k}$ with reversed
orientation) admits, on the complement of a disc, an almost complex structure
$J'$ with
\[ \fo (-M_0^{8k},J')=-\fo (M_0^{8k},J)+\chi (M_0^{8k})\fo (S^{8k})
 = \chi (M_0^{8k})\cdot (1,0).\]
(For the nonstable part $\fo_0$ of the obstruction one can also
deduce this from the formula given in (iii) above.)
Thus, for $a$ and $b$ non-negative integers, the $(4k-1)$-connected manifold
$\#_a M_0^{8k}\#_b (-M_0^{8k})-D^{8k}$ admits an almost complex structure
whose nonstable part of the
obstruction to extension over the closed manifold, by
(i) and (ii), is equal to
\[ \fo_0= a\cdot 0+b\chi (M_0^{8k})-(a+b-1)=b\bigl( \chi (M_0^{8k})-1
\bigr) -(a-1).\]
By a suitable choice of $a$ and $b$, any integer can be realised.
(For spheres of dimension $8k+3$, use the $(4k+1)$-connected 
almost complex $(8k+4)$-manifold constructed in~\cite[Thm.~25]{geig97}.)
The case of $S^7$ is special since ${\mathbb H}P^2-D^8$ admits
an almost complex structure with $\fo =(0,1)\in\bfz\oplus\bfz_2$,
which allows us also to vary the stable part of the obstruction. (The
obstruction class $\fo$ for ${\mathbb H}P^2$ was computed in~\cite{geig97};
instead of the argument employed there one may use the result
of~\cite{heap70} quoted above.) This concludes the proof of
Theorems \ref{thm:contact} and~\ref{thm:contact2}.

\vspace{2mm}

We end this note with a brief discussion of the failure of our
methods to cover stably nontrivial homotopy classes of almost
contact structures on spheres of dimension $\equiv 7$ mod~$8$
(and greater than~$7$). The key to this is the observation that
equation~(\ref{eqn:heaps}) (for $M^{8k}$), which contains purely topological
data, is both necessary and sufficient for {\em any} almost
complex structure on
$M^{8k}-D^{8k}$ to extend as a stable structure over~$M^{8k}$
(this follows from the proof of~\cite[Thm.~1.2]{heap70}).
Thus, in order to construct an exotic contact structure in a stably
nontrivial homotopy class of almost contact structures, our methods
would require us to find a $(4k-1)$-connected $8k$-manifold $M^{8k}$
admitting an almost complex structure on the complement of a disc
and satisfying
\[ \langle \ph (M^{8k})\hat{A}(M^{8k}),[M^{8k}]\rangle\equiv 1
\;\mbox{\rm mod}\; 2. \]
There are two (not entirely unrelated) constructions of such highly
connected manifolds: the Brieskorn construction
described in \cite[Section~5.1]{geig97} and the plumbing of (multiples
of) the tangent disc bundle of $S^{4k}$. Both give rise to {\em parallelisable}
$8k$-manifolds with boundary (for the former this is obvious; for the
latter see~\cite[Satz 8.7]{hima68}). Therefore, Puschnigg's observation shows
that neither of these constructions yields the desired contact
structures. In part (a) of the following proposition we give an
independent proof of this observation for the dimensions of interest
to us. In part (b) we show that
in dimension $15$, at least, the construction via highly
connected almost complex manifolds fails altogether.

\begin{prop}
(a) If $M^{8k}$ is an almost parallelisable
manifold, then the genus
$\langle \ph (M^{8k})\hat{A}(M^{8k}),[M^{8k}]\rangle$
is even. Hence $M^{8k}$ admits a stable almost complex structure.

(b) If $M^{16}$ is a spin manifold whose Pontrjagin numbers other than
$p_2^2$ and $p_4$ vanish, then
$\langle \ph (M^{16})\hat{A}(M^{16}),[M^{16}]\rangle$ is even.
\end{prop}

\begin{proof}
(a) Under the assumptions of the proposition, the only (potentially)
nonvanishing Pontrjagin class of $M^{8k}$ is $p_{2k}(M^{8k})$. Hence,
like in the proof of Lemma~\ref{lem:M_0} we find
\[ \langle \ph (M^{8k})\hat{A}(M^{8k}),[M^{8k}]\rangle =
8k\cdot\hat{A}[M^{8k}]-\frac{p_{2k}[M^{8k}]}{(4k-1)!}.\]
By an argument similar to the proof of the divisibility statement
in Proposition~\ref{prop:acs}, $p_{2k}[M^{8k}]$ is divisible
by $\den (B_{2k}/8k)\cdot (4k-1)!$, cf.~\cite[p.~770]{hirz87}.
Since the denominator of every Bernoulli number is even,
the claim on the parity of the genus follows. As an almost parallelisable
manifold, $M^{8k}$ obviously admits an almost complex structure on
the complement of a disc, which extends as a stable
almost complex structure over the closed manifold
by~\cite[Thm.~1.2]{heap70}. (The condition
$\mbox{\rm Sq}^2 H^{8k-2}(M^{8k})=0$ of that theorem is satisfied
because of the almost parallelisability of~$M^{8k}$.)

(b) By computations analogous to
the proof of Lemma~\ref{lem:M_0} and using the formula (modulo
terms involving $p_1$ or~$p_3$)
\[ \hat{A}(M^{16})=1-\frac{1}{2^5\cdot 3^2\cdot 5}\, p_2(M^{16})+
\frac{1}{2^{16}\cdot 3^4\cdot 5^2\cdot 7}\bigl( 416\, p_2^2(M^{16})
-384\, p_4(M^{16})\bigr) \]
one finds
\[ \langle \ph (M^{16})\hat{A}(M^{16}),[M^{16}]\rangle =
496\, \hat{A}[M^{16}],\]
which is even because of the integrality of the $\hat{A}$-genus
for spin manifolds.
\end{proof}

\begin{ack}
We thank Michael Puschnigg for sending us a copy of his diploma thesis,
and for other interesting correspondence.

F.~D.\ is partially supported by grant no.\ 10201003 of the National Natural
Science Foundation of China. H.~G.\ is partially supported by grant no.\
GE 1245/1-1 of the Deutsche Forschungsgemeinschaft within the
framework of the Schwerpunktprogramm 1154 ``Globale
Differential\-geo\-metrie''.
\end{ack}

\end{document}